\documentstyle[12pt]{article}
\def\R{\hbox{{\rm I}\kern-0.2em{\rm R}\kern0.2em}}%mathematical R for reals
\def\D{\hbox{{\rm I}\kern-0.2em{\rm D}\kern0.2em}}

\def\be{\begin{equation}}
\def\ee{\end{equation}}

\def\({\left(}
\def\){\right)}
\def\[{\left[}
\def\]{\right]}
\def\bc{\begin{center}}
\def\ec{\end{center}}

\begin{document}

{\bf\Large Higher Dimensional Systems of Differential Equations
Obtained by Iterative Use of Complex Methods}

\begin{center}
F M Mahomed$^a$ and Asghar Qadir$^b$\\
$^a$Centre for Differential Equations, Continuum Mechanics and
Applications, School of Computational and Applied Mathematics,
University of the Witwatersrand, Wits 2050, South Africa\\[3ex]
$^b$Centre for Advanced Mathematics and Physics, National
University of Sciences and Technology, Campus H-12, 44000,
Islamabad, Pakistan
\end{center}

{\bf Abstract}.\\

Systems of two ordinary and partial differential equations
(ODEs and PDEs) had been obtained from a scalar complex ODE by
splitting it into its real and imaginary parts. The procedure
was also carried out to obtain a four dimensional system by
splitting a complex system of two ODEs into its real and
imaginary parts. Systems of three ODEs had not been accessible
by these methods. In this paper the complex splitting is used
iteratively to obtain {\it three} and four dimensional systems
of ODEs and four dimensional systems of PDEs for four functions
of {\it two} and four variables. The new systems of four ODEs
are distinct from the class obtained by the single split of a
two dimensional system. Illustrative examples are provided.

\section{Introduction}

Lie \cite{lie} developed his study of the symmetry of
differential equations for complex functions. Of course, to be
differentiable they have to be analytic. Though this was
assumed, analyticity was not used in the sense of treating the
complex function as the real and imaginary parts connected to
each other by the Cauchy-Riemann (CR) conditions. Thus, while
the complex nature of the function was very important for the
topological properties of the Lie groups that arose, it was not
used directly for the differential equations. More recently
\cite{amq1} the splitting of a scalar ODE, into its real and
imaginary parts, was exploited to obtain methods to deal with
systems of ODEs and PDEs. This was called {\it complex symmetry
analysis} (CSA). Of particular interest was its application to
dealing with the variational principle for systems of ODEs
\cite{amq2,faq} and to linearization (conversion of the
equation to linear form by transformation of the dependent and
independent variables) of ODEs \cite{amq3,aqs1,ams}. The latter
was of special interest because of the connection between
geometry and the symmetries of systems of ODEs
\cite{aa,aa1,fmq,aq}. This enabled one to use geometric methods
to linearize systems (including the scalar case) of ODEs
\cite{mq1,mq2}. This method allowed one to not only write down
the linearizing transformation but also to directly provide the
solution. The procedure also led to some new insights regarding
standard systems of equations and to methods for solving
systems that were not amenable to solution by standard symmetry
analysis \cite{faq,aqs2}.

One might have expected that a system of four ODEs or PDEs
could be obtained by using quaternions. This turns out to be
impossible. The reason is that while algebra works for
quaternions calculus does not. Consider $q=w+ix+jy+kz$,
subject to the usual quaternion rules that
$i^2=j^2=k^2=-1~,ij=k=-ji~,jk=i=-kj~,ki=j=-ik$. For $dq/dq=1$
the derivative operator is
\begin{equation}
\frac{d}{dq}=a\frac{\partial}{\partial w}+ib\frac{\partial}{\partial x}
+jc\frac{\partial}{\partial y}+kd\frac{\partial}{\partial z}~,
\end{equation}
subject to the condition that $a-b-c-d=1$. Now ask that
$dq^2/dq=2q$. This requires the above condition along with the
additional conditions, $a-b=a-c=a-d=1$, which are inconsistent
with the earlier requirement. One has to obtain the four
dimensional system by other means, such as using a complex
system of two ODEs and splitting it into a system of four ODEs.

Though one can split a system of two complex ODEs into one of
four ODEs by this procedure, it is not possible to use it to
obtain a system of three ODEs from either a scalar or a vector
equation. In fact, all odd dimensional systems are inaccessible
by the usual CSA splitting procedure. However, one would like to
be able to use the power of CSA for three (and other odd
dimensional systems) dimensional systems as well. In this paper
the complex splitting is used iteratively to be able to generate
three and four dimensional systems of ODEs and to obtain systems
of four PDEs for four functions of two or four variables. The
procedure could be used with more than two iterations to generate
higher dimensional systems as well, but we will not follow that
up here.

The method used to obtain the system of ODEs is to start with a
scalar ODE and regard the dependent variable as a complex
function of a real variable, as was done for CSA, and split
it into a system of two ODEs. Now regard the two dependent
variables as {\it themselves} complex functions of a real
variable. This provides a system of four ODEs. Of course, if
both steps had been merged into one, only a system of two
equations would again have been obtained. (This would be the
reduced system obtainable from the new system of four
equations.) Instead, we first close our eyes to the fact that
we are going to treat the two dependent variables as complex,
and only after obtaining the system of two equations do we
treat each of the dependent variables as themselves complex.
The result of repeating the CSA procedure is different from
treating the dependent variable of the original scalar equation
as a function of two complex functions of a real variable as
the symmetry structure of the systems is
different.

In the CSA procedure it was found that one may be able to
linearize the base scalar equation even though the
corresponding real system does not have enough symmetries to
allow linearization, or even to permit of solution by usual
symmetry methods. The system of four ODEs could then have much
fewer symmetries while the base system is linearizable, or even
have a base system with too few symmetries while the scalar
equation is linearizable \cite{aqs2}. It is to be expected that
here we will not only get the same sort of cases, but that
there would be some even stranger cases arising than for CSA.

One has another option. For the system obtained by the first
split, one can treat one of the two dependent variables as real
and the other as complex. {\it This provides the system of
three ODEs}. In fact, it provides two systems of three ODEs, as
we can choose either of the dependent variables to be real and
the other complex. The two systems obtained are {\it dual} to
each other in some sense. That sense will be clarified by some
examples that we will provide.

Of course, one could take the CSA system of PDEs and treat each
of the dependent variables as complex functions of the two real
variables. This would yield a system of four PDEs for four
functions of two variables. Alternatively, we could have
started with the system of two ODEs and now treated both the
dependent and independent variables as complex. This split also
gives a system of four PDEs for four functions of two
variables. The two systems obtained are dual to each other in a
fairly obvious way. Splitting the system of two PDEs by
treating the dependent and independent variables as complex
gives a system of four PDEs for four functions of four
variables.

The role of the CR-equations was fairly obvious in the original
CSA and had not been spelled out explicitly. A geometric
description of these equations {\it was} given later \cite{aqs2}
but the detailed requirements for a general scalar ODE were not.
On the double use of the splitting procedure the conditions
become more thoroughly coupled and need to be spelled out
explicitly. It turns out that even for the original CSA the
conditions for the derivatives of the functions involved on the
right side of the semilinear equations with respect to the
derivatives of the dependent variables are more complicated
than was envisaged. These have been stated explicitly here.

We have limited the discussion to second order equations only.
One can, of course, go to higher order equations but that
complicates the expressions without providing any further
understanding of the procedures being developed. Also, there is
no direct equivalent of the geometric connection and procedures
for the higher order equations. (There {\it is} an indirect
connection provided by differentiating the second order ODE and
requiring that the original equation hold \cite{mq3,mq4}, but
that will not be followed up here.) On the other hand, one could
have limited oneself to first order equations but that is the
degenerate case and results for it will not hold more generally.

The basic view adopted here is to use the connection of a
second order system of differential equations with a scalar
second order ODE to be able to solve, or otherwise deal with,
the given system. For this purpose, one ideally needs clear
criteria to be able to determine whether the given system is,
or is not, related to some scalar ODE. The criteria should be
such that a computer code could be written to check the given
system for the relationship, and if it is related to construct
the required ODE, which could then be appropriately dealt with.
{\it This had not been done even for the original CSA.} Here we
state the explicit criteria as theorems for the original CSA
and then for the double split systems.

The plan of the paper is as follows. In the next section we
present the basics of the CSA splitting procedure and briefly
mention symmetries of differential equations. We also present
two characterization theorems there that had not been provided
earlier. In the subsequent section we give the split into a
system of three ODEs. In section four we give the split into a
system of four ODEs. The section after that deals with the
split into a system of four PDEs for four functions of two
variables. In section six we present the split into a system of
four PDEs for four functions of four variables. For each of the
systems we provide some examples in the same sections. In the
concluding section we give a brief summary and discussion of the
results.

\section{Complex Splitting and Symmetries}

Consider a general second order ODE
\begin{equation}
u''(r)=f(r;u,u')~.
\end{equation}
We can now take either $u$ to be a complex function of the real
variable $r$ or take both $u$ {\it and} $r$ to be complex. Let
$u=p+\iota q$ in the former case and put
\begin{equation}
f(r;u,u')=f^r(r;p,q;p',q')+\iota f^i(r;p,q;p',q')~.
\end{equation}
The resulting system of ODEs is
\begin{equation}
p''(r)=f^r(r;p,q;p',q')~,q''(r)=f^i(r;p,q;p',q')~.
\end{equation}
The CR-equations for this system are
\begin{equation}
f^r_p=f^i_q~,f^r_q=-f^i_p~;
f^r_{p'}=f^i_{q'}~;f^r_{q'}=-f^i_{p'}~;
\end{equation}
with no conditions on $p$ and $q$ other than second
differentiability.

For the latter case let $r=s+\iota t$ as well. Then
\begin{equation}
f(r;u,u')=f^r(s,t;p,q;p_s,q_s,p_t,q_t)+
\iota f^i(s,t;p,q;p_s,q_s,p_t,q_t)~.
\end{equation}
Splitting the equation into its real and imaginary parts then
gives the system of two PDEs
\begin{equation}
p_{ss}-p_{tt}+2q_{st}=4f^r~,~ q_{ss}-q_{tt}-2p_{st}=4f^i~.
\end{equation}
The CR-equations for this system include conditions for $p$ and
$q$, apart from the previous ones (which are now more
complicated),
\begin{eqnarray}
p_s=q_t~,~p_t=-q_s~; \nonumber\\
f^r_s=f^i_t~,~f^r_t=-f^i_s~; \nonumber\\
f^r_p=f^i_q~,~f^r_q=-f^i_p~.
\end{eqnarray}
The derivative of the functions with respect to the derivatives
is more complicated to state. The problem is that
\begin{equation}
2u'\rightarrow (\frac{\partial}{\partial s}
-\iota \frac{\partial}{\partial t})(p+\iota ~q)
=(p_s+q_t)-\iota (p_t-q_s)~.
\end{equation}
Thus, for $f(r;u,u')$ to be analytic, $f^r$ and $f^i$ cannot
depend arbitrarily on $p_s$, $q_s$, $p_t$ and $q_t$ but must
depend on $p_s+q_t$ and $p_t-q_s$. If we call these
variables $\phi$ and $\psi$, respectively, then the last
conditions are
\begin{eqnarray}
f^r_{\phi}=f^i_{\psi}~,~f^r_{\psi}=-f^i_{\phi}~.
\end{eqnarray}

Though the CR-equations were taken as obvious, they are needed
to characterize systems of (real) ODEs and PDEs that can arise
by splitting a scalar ODE by CSA methods. To complete the CSA
procedure we state the following two theorems.

{\bf Theorem 1:} {\it A system of two second order ODEs} (4)
{\it corresponds to a scalar second order ODE} (2) {\it if and
only if it satisfies the system of CR-conditions} (5).

{\bf Theorem 2:} {\it A system of two second order PDEs} (7)
{\it corresponds to a scalar second order ODE} (2) {\it if and
only if it satisfies the system of CR-conditions} (8) {\it and}
(10) {\it where $\phi=p_s+q_t$ and $\psi=p_t-q_s$.}

Symmetry analysis deals with the infinitesimal generators that
leave the differential equation invariant under point
transformations $(r,u)\rightarrow(R,U)$ say,
\begin{equation}
{\bf X}=\xi(r,u)\frac{\partial}{\partial r}+\eta(r,u)\frac{\partial}{\partial u}~.
\end{equation}
To be able to apply the operators to differential equations
these generators have to be {\it prolonged} or {\it extended}
to include the higher derivatives. For the first derivative
\begin{equation}
{\bf X}^{[1]}=\xi \frac{\partial}{\partial r}+\eta \frac{\partial}{\partial u}
+\eta^{[1]}(r;u,u')\frac{\partial}{\partial u'},
\end{equation}
where
\begin{equation}
\eta^{[n]}=\frac{d\eta^{[n-1]}}{dr}-u^{n}(r)\frac{d\eta}{dr}~,
\end{equation}
($\eta^{[0]}=\eta$), $d/dr$ stands for the total
derivative,
\begin{equation}
\frac{d}{dr}=\frac{\partial}{\partial r}+u'\frac{\partial}
{\partial u}+u''\frac{\partial}{\partial u'}
\end{equation}
and $u^{n}(r)$ stands for the $n^{th}$ derivative. The
prolonged generators for higher order differential equations
can be similarly obtained using (13).

We can extend the analysis to systems of ODEs by replacing the
scalar $u$ by a vector ${\bf u}$ and the corresponding partial
derivative by $\nabla_{\bf u}$. Consequently we must replace
the scalar $\eta$ and $\eta^{[1]}$ by the vectors
$\underline{\eta}$ and $\underline{\eta}^{[1]}$. The extension
to the PDEs is more complicated. The scalar $r$ has to now also
be replaced by a vector $\bf s$ and along with it the scalar
$\xi$ by the vector $\underline{\xi}$, but now the derivative
of $\bf u$ becomes $\nabla_{\bf s}\bf u$ and the derivative with
respect to this vector of partial derivatives becomes too messy
to read. As such, we write $\nabla_{\bf s}{\bf u}={\bf u_{1}}$.
Of course, we need to also bear in mind that the variables that
the functions depend on will not be the derivatives given but linear
combinations as we saw for the CR-conditions. The second
derivative can then be written as ${\bf u_{2}}$ and so on. For
the complex case the ${\bf X}$ was written as ${\bf Z}$ and
similar changes were made for the coefficients but it will be
more convenient to use the same notation throughout here.

\section{System of Three ODEs by Double Splitting}

Consider (2) as the base scalar equation with the split of the
function given by (3) leading to the system (4) subject to the
CR-equations (5). Now regard $p$ as a real variable $x$ and $q$
as the complex variable $y+\iota z$. We run into a problem
here. There are three second order ODEs but the number of
functions to be obtained from $f$ must be even. To avoid this
problem we put
\begin{equation}
f(r;u,u')=g(r;p,p')+\iota ~G(r;p,q;p,q')~.
\end{equation}
The choice of what part to put into $g$ and what part to put
into $G$ is clearly arbitrary. For definiteness, we define $g$
to consist of {\it all} those terms that do not involve $q$ or
$q'$. Then $G$ consists of all those terms that do. Now we
proceed with the second split by putting
\begin{equation}
G(r;p,q;p,q')=k(r;x,y,z;x',y',z')+\iota ~l(r;x,y,z;x',y',z')~,
\end{equation}
yielding the system of three ODEs
\begin{equation}
x''=h(r;x,x')~,~y''=k(r;{\bf x},{\bf x'})~,z''=l(r;{\bf x},{\bf x'})~,
\end{equation}
subject to the CR-equations
\begin{equation}
k_y=l_z~,~k_z=-l_y~;~k_{y'}=l_{z'}~,~k_{z'}=-l_{y'}~.
\end{equation}
The first prolonged symmetry generator is:
\begin{eqnarray}
{\bf X}^{[1]}=\xi(r,{\bf x})\frac{\partial}{\partial r}+\eta^x(r,{\bf x})
\frac{\partial}{\partial x}+\eta^y(r,{\bf x})\frac{\partial}{\partial y}
+\eta^z(r,{\bf x})\frac{\partial}{\partial z} \nonumber\\
+\eta^{x[1]}(r;{\bf x},\nabla_{{\bf x}'})\frac{\partial}{\partial x'}
+\eta^{y[1]}(r;{\bf x},\nabla_{{\bf x}'})\frac{\partial}{\partial y'}
+\eta^{z[1]}(r;{\bf x},\nabla_{{\bf x}'})\frac{\partial}{\partial z'}~.
\end{eqnarray}

Instead of this procedure, at the second split we could have
taken $p(r)=x(r)+\iota ~y(r)$ and set $q(r)=z(r)$. This will
give a {\it dual system} in some sense. We would now have to
set
\begin{equation}
f(r;u,u')=g(r;x,y,z;x',y',z')+\iota ~h(r;x,y,z;x',y',z')+\iota ~k(r;z,z')~,
\end{equation}
and the slightly modified CR-equations
\begin{equation}
g_y=h_z~,~g_z=-h_y~;~g_{y'}=h_{z'}~,~g_{z'}=-h_{y'}~.
\end{equation}
The prolonged symmetry generator remains unaltered in form except
that now the coefficient of the last term is pure imaginary, to
account for the last imaginary term in (20). The sense of the
duality will be clarified by the examples.

For completeness we state a theorem for the characterization of
systems of three ODEs that correspond to a scalar ODE by double
splitting.

{\bf Theorem 3:} {\it The system of three ODEs} (17) {\it
corresponds to the scalar ODE} (2){\it , for any consistent
identification of the function given by} (15) {\it and}
(16){\it , provided the CR-conditions} (18) {\it hold.}

{\bf Remark:} {\it The dual procedure gives the same system but
now we require that the split given by} (20) {\it and} (21)
{\it holds.}

{\bf Example 1:} Consider the free-particle scalar ODE, $u''=0$.
It clearly yields the system $x''=y''=z''=0$. The splitting of
the functions is obviously trivial. However, the infinitesimal
symmetry generators of the system are not trivially related to
the generators of the scalar ODE. Even for the original CSA, it
had been noted that the symmetries for the system could not be
a simple doubling of the symmetries of the original ODE, as the
maximal Lie algebra for the system is $sl(4,\R)$, which has 15
generators, while the algebra for the scalar free particle
equation is $sl(3,\R)$, which has 8. Doubling gives one extra
generator. It is also clear that it cannot be a simple matter of
leaving one generator out, as the complex generators occur in
pairs. What happened there was that we lost two generators and
got one new one. The system must clearly have 24 generators, as
the algebra is $sl(5,\R)$. However, the double splitting gives
only 23 operators of which 15 are symmetries and 8 are Lie-like
\cite{amq1}, five of the nine dilations coming from the dependent
variable and the four local projective symmetries are missing.
We have to use the actual symmetries obtained here and require
closure of the algebra to generate the full 24. The dual system
also gives the same symmetry structure, as expected. However,
due to the fact that the $k$ in (20) has a coefficient with
iota, {\it even though $k$ itself is zero}, the operator for
the symmetry carries an imaginary. If we do not put in the
iota for the operator, we lose some of the symmetry generators.

{\bf Example 2:} Consider the scalar equation \cite{pl}
\begin{equation}
u''=s^{-5}u^2 ~.
\end{equation}
It splits into the following system of ODEs:
\begin{eqnarray}
x''=-2s^{-5}yz~; \nonumber \\
y''=-2s^{-5}zx~; \nonumber \\
z''=-2s^{-5}xy~;
\end{eqnarray}
subject to the further algebraic constraint
\begin{equation}
x^2+y^2=z^2~.
\end{equation}
There are no symmetries among the 8 Lie-like operators
obtained from the splitting of the complex generators
\begin{equation}
{\bf X}_1=s\frac{\partial}{\partial s}+3u\frac{\partial}{\partial u}~;~
{\bf X}_2=s^2\frac{\partial}{\partial s}+su\frac{\partial}{\partial u}~.
\end{equation}
However, the system admits the two symmetry generators:
\begin{eqnarray}
{\bf Y}_1=x\frac{\partial}{\partial x}+y\frac{\partial}{\partial y}
+z\frac{\partial}{\partial z}~; \nonumber \\
{\bf Y}_2=s^2\frac{\partial}{\partial s}+sx\frac{\partial}{\partial x}
+sy\frac{\partial}{\partial y}+sz\frac{\partial}{\partial z}~.
\end{eqnarray}

{\bf Example 3:} Consider the scalar (Emden-Fowler)
equation \cite{pr}
\begin{equation}
u''+5s^{-1}u'+u^2=0~.
\end{equation}
It splits into the following system of ODEs:
\begin{eqnarray}
x''+5s^{-1}x'-2yz=0~; \nonumber \\
y''+5s^{-1}y'-2zx=0~; \nonumber \\
z''+5s^{-1}z'-2xy=0~;
\end{eqnarray}
subject to the further algebraic constraint
\begin{equation}
x^2+y^2=z^2~.
\end{equation}
The Emden-Fowler equation has only the one scaling
symmetry
\begin{equation}
{\bf X}=s\frac{\partial}{\partial s}-2u\frac{\partial}{\partial u}~.~
\end{equation}
The split system has 4 Lie-like operators, none of which are
symmetries of the system. However, the system admits the scaling
\begin{equation}
{\bf Y}=s\frac{\partial}{\partial s}-2x\frac{\partial}{\partial x}-2y\frac{\partial}{\partial y}
-2z\frac{\partial}{\partial z}~.
\end{equation}

\section{System of Four ODEs by Double Splitting}

For the four dimensional system, after the first split of (2)
given by (3), (4) and (5) we could set $p(r)=w(r)+\iota ~x(r)$
and $q(r)=y(r)+\iota ~z(r)$ and
\begin{eqnarray}
f^r(r;p,q;p',q')=g(r;{\bf w},{\bf w}')+\iota ~h(r;{\bf w},{\bf w}')~, \nonumber\\
f^i(r;p,q;p',q')=k(r;{\bf w},{\bf w}')+\iota ~l(r;{\bf w},{\bf w}')~,
\end{eqnarray}
yielding the system of four ODEs
\begin{eqnarray}
w''(r)=g(r;{\bf w},{\bf w}')~,~x''(r)=h(r;{\bf w},{\bf w}')~, \nonumber\\
y''(r)=k(r;{\bf w},{\bf w}')~z''(r)=l(r;{\bf w},{\bf w}')~,
\end{eqnarray}
subject to the CR-conditions
\begin{eqnarray}
g_w+h_x=k_y+l_z~,~g_x-h_w=k_z-l_y~, \nonumber\\
g_y+h_z=-k_w-l_x~,~g_z-h_y=-k_x+l_w~, \nonumber\\
g_{w'}+h_{x'}=k_{y'}+l_{z'}~,~g_{x'}-h_{w'}=k_{z'}-l_{y'}~, \nonumber\\
g_{y'}+h_{z'}=-k_{w'}-l_{x'}~,~g_{z'}-h_{y'}=-k_{x'}+l_{w'}~,
\end{eqnarray}
where ${\bf w}=(w,x,y,z)$. The prolonged symmetry generator can
now be written as
\begin{equation}
{\bf X}=\xi(r,{\bf w})\frac{\partial}{\partial r}+\underline{\eta}
(r,{\bf w}){\bf .}\nabla_{\bf w}+\underline{\eta}^{[1]}
(r;{\bf w},{\bf w}'){\bf .}\nabla_{{\bf w}'}~.
\end{equation}
Writing this equation out in detail makes it too unwieldy to
convey much wisdom.

There is no dual system to this as the splitting is direct. We
could have obtained a system of four ODEs by a three stage
splitting, setting one of the dependent variables in the
systems of three ODEs as a complex pair. For each of the three
dimensional systems obtained, one splitting would give the
above system and one would give a new system. We are, here,
limiting our discussion to a two-step splitting only.

The characterization theorem is: \\
{\bf Theorem 4:} {\it The system of four second order ODEs}
(33) {\it corresponds to the scalar ODE} (2) {\it by double
complex splitting if and only if the CR-conditions} (34)
{\it hold with the splitting} (32).

The examples will illustrate our procedure further.

{\bf Example 4:} The free particle scalar equation obviously
yields the free-particle system of four equations,
$w''=x''=y''=z''=0$. The symmetry algebra must be $sl(6,\R)$,
which has 35 generators. There are a total of 18 symmetries and
8 Lie-like operators. Again, the missing ones come from the
dilations involving the dependent variables and local projective
symmetries. Again, the closure of the algebra starting with the
derived symmetries generates the full $sl(6,\R)$.

{\bf Example 5:} Consider (22) and now split into the system of
four ODEs:
\begin{eqnarray}
w''=s^{-5}(w^2-x^2-y^2+z^2)~; \nonumber \\
x''=s^{-5}(2wx- 2yz)~; \nonumber \\
y''=s^{-5}(2wy-2xz)~; \nonumber \\
z''=s^{-5}(2wz+2xy)~.
\end{eqnarray}
This has 8 Lie-like operators of which none are symmetries of
the system. However, the system does admit the two symmetries
\begin{eqnarray}
{\bf Y}_1=s\frac{\partial}{\partial s}+3w\frac{\partial}{\partial w}
+3x\frac{\partial}{\partial x}+3y\frac{\partial}{\partial y}
+3z\frac{\partial}{\partial z}~; \nonumber \\
{\bf Y}_2=s^2\frac{\partial}{\partial s}+sw\frac{\partial}{\partial w}+sx\frac{\partial}{\partial x}
+sy\frac{\partial}{\partial y}+sz\frac{\partial}{\partial z}~.
\end{eqnarray}

{\bf Example 6:} Again consider the scalar (Emden-Fowler)
equation (27). It splits into the following system of ODEs:
\begin{eqnarray}
w''+5s^{-1}w'+w^2-x^2-y^2+z^2=0~; \nonumber \\
x''+5s^{-1}x'+2wx-2yz=0~; \nonumber \\
y''+5s^{-1}y'+2wy-2xz=0~; \nonumber \\
z''+5s^{-1}z'+2wz+2xy=0~.
\end{eqnarray}
The split system again has 4 Lie-like operators, none of which
are symmetries of the system. However, the system admits the
scaling
\begin{equation}
{\bf Y}=s\frac{\partial}{\partial s}-2w\frac{\partial}{\partial w}
-2x\frac{\partial}{\partial x}-2y\frac{\partial}{\partial y}
-2z\frac{\partial}{\partial z}~.
\end{equation}

\section{System of Four PDEs for Four Functions of Two
Variables by Double Splitting}

For the PDEs, we need to treat the dependent variable as
complex. However, there are two ways of doing so if we are to
obtain functions of only two variables. We could retain the
independent variable as real in the first step and then make it
complex in the second step, or first treat it as complex and
then retain the same independent variables in the second step.
We follow the former procedure first. Starting with (2), with
the first step splitting given by (3), (4) and (5), we put
$r=s+\iota ~t$ and proceed with treating the dependent
variables as complex exactly as in the case for the system of
four ODEs.

The system of equations is now
\begin{eqnarray}
w_{ss}-w_{tt}+2x_{st}=4g({\bf s};{\bf w};{\bf w}_s,{\bf w}_t)~,
~x_{ss}-x_{tt}-2w_{st}=4h({\bf s};{\bf w};{\bf w}_s,{\bf w}_t)~; \nonumber\\
y_{ss}-y_{tt}+2z_{st}=4k({\bf s};{\bf w};{\bf w}_s,{\bf w}_t)~,
~z_{ss}-z_{tt}-2y_{st}=4l({\bf s};{\bf w};{\bf w}_s,{\bf w}_t)~,
\end{eqnarray}
subject to the CR-equations
\begin{eqnarray}
w_s=x_t~,~w_t=-x_s~,~y_s=z_t~,~y_t=-z_s~; \nonumber\\
g_w=h_x~,~g_x=-h_w~,~g_y=h_z~,g_z=-h_y~; \nonumber\\
k_w=l_x~,k_x=-l_w~,~k_y=l_z~,~k_z=-l_w~.
\end{eqnarray}

Now, as for (10), defining
\begin{equation}
\phi=w_s+x_t~,~\psi=w_t-s_x~,~\kappa=y_s+z_t~,~\lambda=y_t-z_s~,
\end{equation}
the conditions for the derivatives with respect to the
derivatives can be written as
\begin{equation}
g_{\phi}=h_{\psi}~,~g_{\psi}=-h_{\phi}~;~k_{\kappa}=l_{\lambda}~,
~k_{\lambda}=-l_{\kappa}~.
\end{equation}

Writing $(s,t)=\underline{s}$, for the infinitesimal symmetry
generator we will now have two components for
$\underline{\xi}$, namely $(\xi^s,\xi^t)$ and
$\underline{\eta}$ will have four components as for the system
of four ODEs. The additional feature is that the prolonged
derivatives will be for $\nabla_{{\bf w}_{\bf s}}$ and the
coefficients will be $\underline{\eta}^{{\bf [1]}}$. Thus
the prolonged generator can be written as
\begin{equation}
{\bf X}=\underline{\xi}({\bf s},{\bf w}){\bf .}\nabla_{\bf s}+
\underline{\eta}({\bf s},{\bf w}){\bf .}{\nabla_{\bf w}}+
\underline{\eta}^{{\bf [1]}}({\bf s},{\bf w})
{\bf .}\nabla_{\nabla_{\bf s}{\bf w}}~,
\end{equation}
where $\underline{\eta}^{\nabla_{{\bf [1]}}}({\bf s},{\bf w})$
is the generalization of $\underline{\eta}^{[1]}$ for the
case of PDEs.

The characterization theorem here is: \\
{\bf Theorem 5:} {\it The system of four PDEs for four
functions of two variables} (40) {\it corresponds to the scalar
ODE} (2) {\it by double complex splitting, if and only if the
CR-conditions} (41) - (43) {\it hold, provided that $g,h,k,l$
depend on the derivatives only in the combinations given by}
(43).

Instead, if we had put $r=s+\iota ~t$ in the first step and
then proceeded to split the dependent variables twice, we would
have got the dual system
\begin{eqnarray}
w_{ss}-w_{tt}+2y_{st}=4g({\bf s};{\bf w};{\bf w}_x,{\bf w}_t)~,
~x_{ss}-x_{tt}+2z_{st}=4h({\bf s};{\bf w};{\bf w}_x,{\bf w}_t)~; \nonumber\\
y_{ss}-y_{tt}-2w_{st}=4k({\bf s};{\bf w};{\bf w}_x,{\bf w}_t)~,
~z_{ss}-z_{tt}-2x_{st}=4l({\bf s};{\bf w};{\bf w}_x,{\bf w}_t)~.
\end{eqnarray}
The CR-conditions are considerably more involved. The simple
ones are
\begin{eqnarray}
w_s=y_t~,~w_t=-y_s~,~x_s=z_t~,~x_t=-z_s~; \nonumber\\
g_w+h_x=k_y+l_z~,~g_x-h_w=k_z-l_y~; \nonumber\\
g_y+h_z=-k_w-l_x~,g_z-h_y=-k_x+l_w~.
\end{eqnarray}
For the derivatives with respect to derivatives we have to
define the new variables
\begin{equation}
\alpha =w_s+y_t~,~\beta=x_s+z_t~,~\gamma=w_t-y_s~,~\delta=x_t-z_s~
\end{equation}
to get
\begin{eqnarray}
g_{\alpha}+h_{\beta}=k_{\gamma}+l_{\delta}~,~g_{\beta}-h_{\alpha}
=k_{\delta}-l_{\gamma}~, \nonumber\\
g_{\gamma}+h_{\delta}=-k_{\alpha}-l_{\beta}~,~g_{\delta}-h_{\gamma}
=-k_{\beta}+l_{\alpha}~.
\end{eqnarray}
The form of the symmetry generator remains unchanged. Once
again, we rely on the examples to illustrate the procedure.

Here we need a separate theorem because the systems are
apparently different, though they are dual to each other in
some sense. \\
{\bf Theorem 6:} {\it The} ({\it dual}) {\it system of four
PDEs for four functions of two variables} (45) {\it corresponds
to the scalar ODE} (2) {\it by double complex splitting, if and
only if the CR-conditions} (46) - (48) {\it hold, provided that
$g, h, k, l$ only depend on the derivatives in the combinations
given by} (48).

{\bf Example 7:} Consider the free particle equation split into
the system of four PDEs for two independent variables
\begin{eqnarray}
w_{ss}-w_{tt}+2x_{st}=0~; \nonumber \\
x_{ss}-x_{tt}-2w_{st}=0~; \nonumber \\
y_{ss}-y_{tt}+2z_{st}=0~; \nonumber \\
z_{ss}-z_{tt}-2y_{st}=0~.
\end{eqnarray}
The symmetry generators split into 28 Lie-like operators of
which 20 are symmetries. However, the system admits infinitely
many symmetries.

The dual system is very similar. More precisely, it is the the
system (45) with the right side set equal to zero. Since the
CR-conditions are trivial there is no significant difference
between the original and the dual system.

{\bf Example 8:} Again consider the scalar (Emden-Fowler)
equation (27). It splits into the following system of four PDEs of
two independent variables:
\begin{eqnarray}
w_{ss}-w_{tt}+2x_{st}+10\frac{s}{s^2+t^2}(w_s+x_t)+10\frac{t}{s^2+t^2}(x_s-w_t) \nonumber \\
+4(w^2-x^2-y^2+z^2)=0~; \nonumber \\
x_{ss}-x_{tt}-2w_{st}+10\frac{s}{s^2+t^2}(x_s-w_t)
-10\frac{t}{s^2+t^2}(w_s+x_t)+8(wx-yz)=0~; \nonumber \\
y_{ss}-y_{tt}+2z_{st}+10\frac{s}{s^2+t^2}(y_s+z_t)
+10\frac{t}{s^2+t^2}(z_s-y_t)+8(wy-xz)=0~; \nonumber \\
z_{ss}-z_{tt}-2y_{st}+10\frac{s}{s^2+t^2}(z_s-y_t)
-10\frac{t}{s^2+t^2}(y_s+z_t)+8(wz+xy)=0~.\end{eqnarray}
This split system again has 4 Lie-like operators, none of which are
symmetries of the system. However, it admits the scaling
\begin{equation}
{\bf Y}=s\frac{\partial}{\partial s}+t\frac{\partial}{\partial t}
-2w\frac{\partial}{\partial w}-2x\frac{\partial}{\partial x}
-2y\frac{\partial}{\partial y}-2z\frac{\partial}{\partial z}~.
\end{equation}

\section{System of Four PDEs for Four Functions of Four
Variables by Double Splitting}

This is the most straight forward (and the most complicated) of
the various possibilities considered. At the first step we
regard both the independent and the dependent variables as given
by (6) to (8). For the second step we run short of symbols for
the variables. As such, we now write the first variable
(previously written as $s$), as the complex variable $s+\iota ~t$
and the second variable (previously written as $t$), as the
complex variable $u+\iota ~v$ and write $\bf s$ for $(s,t,u,v)$.
Further, we put
\begin{equation}
p(s,t)\rightarrow w({\bf s})+\iota ~x({\bf s})~,
~q(s,t)\rightarrow y({\bf s})+\iota ~z({\bf s})~,
\end{equation}
\begin{eqnarray}
f^r(s,t;p,q;p_s,q_s,p_t,q_t)=g({\bf s};{\bf w},\nabla_{{\bf s}}{\bf w})
+\iota ~h({\bf s};{\bf w},\nabla_{{\bf s}}{\bf w})~; \nonumber\\
f^i(s,t;p,q;p_s,q_s,p_t,q_t)=k({\bf s};{\bf w},\nabla_{{\bf s}}{\bf w})
+\iota ~l({\bf s};{\bf w},\nabla_{{\bf s}}{\bf w})~.
\end{eqnarray}

The system of equations is
\begin{eqnarray}
w_{ss}-w_{tt}+2x_{st}-w_{uu}+w_{vv}-2x_{uv}+2y_{su}-2y_{tv} \nonumber\\
+2z_{sv}+2z_{tv}=4g({\bf s};{\bf w},\nabla_{{\bf s}}{\bf w}); \nonumber\\
x_{ss}-x_{tt}-2w_{st}-x_{uu}+x_{vv}+2w_{uv}+2z_{su}-2z_{tv} \nonumber\\
-2y_{sv}-2y_{tv}=4h({\bf s};{\bf w},\nabla_{{\bf s}}{\bf w}); \nonumber\\
y_{ss}-y_{tt}+2z_{st}-y_{uu}+y_{vv}-2z_{uv}+2w_{su}-2w_{tv} \nonumber\\
+2x_{sv}+2x_{tv}=4k({\bf s};{\bf w},\nabla_{{\bf s}}{\bf w}); \nonumber\\
z_{ss}-z_{tt}-2y_{st}-z_{uu}+z_{vv}+2y_{uv}+2x_{su}-2x_{tv} \nonumber\\
-2w_{sv}-2w_{tv}=4l({\bf s};{\bf w},\nabla_{{\bf s}}{\bf w});~
\end{eqnarray}
subject to the CR-conditions
\begin{eqnarray}
w_s+x_t=y_u+z_v~,~w_t-x_s=y_v-z_u~, \nonumber\\
w_u+x_v=-y_s-z_t~,~w_v-x_u=-y_v+z_u~; \nonumber\\
g_s+h_t=k_u+l_v~,~g_t-h_s=k_v-l_u~, \nonumber\\
g_u+h_v=-k_s-l_t~,~g_v-h_u=-k_t+l_s~; \nonumber\\
g_w+h_x=k_y+l_z~,~g_x-h_w=k_z-l_y~, \nonumber\\
g_y+h_z=-k_w-l_x~,~g_z-h_y=-k_x+l_w~.
\end{eqnarray}
The derivatives with respect to the derivatives require the
variables
\begin{eqnarray}
\alpha=w_s+x_t+y_u+z_v~,~\beta=w_t-x_s+y_v-z_u~; \nonumber\\
\gamma=w_u+x_v-y_s-z_t~,~\delta=w_v-x_u-y_t+z_s~.
\end{eqnarray}
Then the rest of the CR-conditions are
\begin{eqnarray}
g_{\alpha}-h_{\beta}=k_{\gamma}-l_{\delta}~,~g_{\beta}+
h_{\alpha}=k_{\delta}+l_{\gamma}~; \nonumber\\
g_{\gamma}-h_{\delta}=-k_{\alpha}+l_{\beta}~,~g_{\delta}+
h_{\gamma}=-k_{\beta}-l_{\alpha}~.
\end{eqnarray}

The characterization theorem in this case is: \\
{\bf Theorem 7:} {\it The system of four second order PDEs
for four functions of four variables} (54) {\it corresponds
to the scalar second order ODE} (2) {\it by double complex
splitting provided the functions $g, h, k, l$ depend on the
derivatives only in the combinations given by} (57) {\it
and the CR-conditions} (55) {\it and} (57) {\it hold.}

The prolonged symmetry generator for the system is
\begin{eqnarray}
{\bf X}^{[1]}=\underline{\xi}({\bf s},{\bf g}){\bf .}\nabla_{{\bf s}}
+\underline{\eta}({\bf s},{\bf g}){\bf .}\nabla_{{\bf g}}+
\underline{\eta}^{{\bf [1]}}({\bf s},{\bf g},\nabla_{{\bf s}}{\bf g){\bf .}
\nabla_{\nabla_{{\bf s}}}{\bf g}}~.
\end{eqnarray}
The derivatives with respect to derivatives are to be taken
bearing in mind the discussion for the CR-equations.
However, even if we ignore it in taking the derivatives, no
error will ensue.

We again rely on the examples to illustrate our systems.

{\bf Example 9:} The free-particle system of equations is given
by (54), with the right side set equal to zero. The
CR-conditions are trivial. There are now 32 Lie-like operators,
of which only 24 are symmetry generators. As before the local
projective symmetries are lost. Here there are 8 such. However,
the dilations are {\it not} lost here. The system, itself, has
an infinite number of symmetry generators.

{\bf Example 10:} Consider the double splitting of the
Emden-Fowler equation (27) into a system of four PDEs for four
functions of four variables. It is given by (54) with
\begin{eqnarray}
g=C[(sA+tB)\alpha+(tA-sB)\beta+(uA+vB)\gamma+(vA-uB)\delta] \nonumber \\
-w^2+x^2+y^2-z^2~; \nonumber \\
h=C[(sA+tB)\beta-(tA-sB)\alpha+(uA+vB)\delta-(vA-uB)\gamma] \nonumber \\
-2wx+2yz~; \nonumber \\
k=C[(sA+tB)\gamma+(tA-sB)\delta-(uA+vB)\alpha-(vA-uB)\beta] \nonumber \\
-2wy+2xz~; \nonumber \\
l=C[-(sA+tB)\delta+(tA-sB)\gamma+(uA+vB)\beta-(vA-uB)\alpha] \nonumber \\
-2wz-2xy~;
\end{eqnarray}
where
\begin{equation}
C=-\frac{5}{A^2+B^2}~,~A=s^2-t^2+u^2-v^2~,~B=2st+2uv~.
\end{equation}
There are 4 Lie-like operators none of which are symmetries.
The system, like the scalar equation, has a scaling symmetry.

\section{Summary and Discussion}

In this paper we have considered systems of three and four
second order ODEs, and systems of four second order PDEs for
four functions of two or four variables, that correspond to a
scalar equation, that we shall call a {\it base equation} by a
specific procedure, that we call {\it double complex splitting}.
We have also provided characterization criteria for such
systems to correspond to the base equation and a clear procedure
to be able to construct the base equation. Thus, in principle,
we could write a computer code that could take any such system
given and check if it corresponds to a base equation. It could
then construct the base equation.

What is the advantage of having such base equations and
constructing them? The point is that it is much easier to deal
with the base equation than with the system. Thus, for example,
if the base system has two infinitesimal symmetry generators it
could be solved by symmetry methods. In fact it could have
eight symmetry generators and thus be linearizable. In that
case we could write down the solution directly. Following the
double splitting procedure by which the system corresponds, we
could then write down the solution for the system of ODEs or
PDEs. Note that in this procedure the system need not have the
required symmetry for being directly solvable by symmetry
methods.

The procedure adopted for CSA could only give an even
dimensional system as it simply split $n$ (complex) equations
into $2n$ (real) equations. Also, the number of independent
variables in the PDE equals the number of dependent variables.
In the complex double split procedure we are being considerably
more adventurous. Having obtained the system of two ``real"
equations we conveniently forget that they arose from a scalar
complex equation, treat it as complex and promptly split the
equations again. Now we have the earlier freedom of choosing the
independent variable to be either real or complex, while
treating the dependent variables to be complex, {\it but we have
the additional freedom to choose to treat one of the dependent
variables as real and the other as complex}. This provides the
possibility of obtaining an odd dimensional system. Further, to
obtain the PDEs, we could choose to treat the independent
variables as real first and then complex, complex first and then
real or complex both times. Thus we also get the system of four
PDEs for four functions of two variables.

In the cases of full double complex splitting, where either the
dependent variables were fully split or both the dependent and
independent variables were double split, giving the system of
four ODEs or four PDEs for four functions of four variables,
there were no complications of additional dual systems arising.
However, in the case of the system of three ODEs or the system
of four PDEs for four functions of two variables, we got dual
systems arising. The duality in the former case was very
obvious but in the latter it was considerably more involved on
account of the CR-conditions. Even in the case of the system of
three ODEs the symmetry generators had to carry an iota to make
sense. The ``duality" of these systems needs to be better
understood. Note that $\alpha$ in (56) is simply $\nabla_{{\bf
s}}{\bf .}{\bf w}$. It would be interesting to find out what
the operators for the other variables are. Presumably, they
would be ``dual" divergence operators in some sense. This may
shed some light on the structure of the double split systems.

For the system of three ODEs obtained by double splitting, it
would be of interest to consider the ambiguity due to the
choice of the function $g$ in (15), or $k$ in (20) for the dual
system. In some sense all choices must be ``equivalent". The
question is whether one gets an equivalence class. Further,
{\it would they be equivalent under point transformations} or
possibly some more general transformations like contact or
higher order transformations \cite{ni}.

The algebraic constraint that arises in the system of three
ODEs was not apparent in setting up the system but was found
in the examples. It is interesting to note that it
geometrically amounts to the solution lying on a cone. This
seems to be generic for the three dimensional system. It also
shows that it can be written as a system of two ODEs. However,
that system is much more complicated.

One would have hoped that for the system of PDEs corresponding
to a base ODE one could use the symmetries of the ODE to
obtain a ``core" set of symmetries for the system of PDEs.
Even if the system has infinitely many symmetries, the base
equation can only have a finite number. However, the examples
show that we can lose all the Lie symmetry generators and be
left only with Lie-like ones. This applies even to the free
particle equation. Further for the Emden-Fowler equation, we
are left with {\it no} Lie symmetries from the Lie-like
operators, though the equation has a scaling symmetry and so
does the double-split system. In general, we obtain {\it
Lie-like operators} and not {\it Lie-symmetry generators} that
would form an algebra. The Lie-like operators somehow encode
the symmetries of the base equation. It would be most
important to learn {\it how} they do so. It may be that the
CR-conditions will enable us to re-construct the Lie from the
Lie-like symmetry.

It is of interest to note that not only for the PDEs but also
for the systems of ODEs, we get Lie-like operators arising and
lose some Lie-symmetry generators. It would be worth while to
see this encoding of symmetry as distinct from the PDE case.

We hope that in future the use of these systems for the
variational principle and with linearization will be
followed up. It should lead to interesting and useful
results.

\section*{Acknowledgments}
AQ is most grateful to DECMA and CAM of Wits for support during
a visit in 2011.

\end{document}